\documentclass[12pt]{article}
\def\finproof{\hfill\hbox{\vrule width1.0ex height1.5ex}\vspace{2mm}}

\begin{document}

\begin{center}
{\LARGE 
Non-linear partial differential equations with discrete
state-dependent delays in a metric space }

\bigskip \bigskip

{\sc Alexander V. Rezounenko}

\smallskip \bigskip

Department of Mechanics and Mathematics, Kharkov
University,

 4, Svobody Sqr., Kharkov, 61077, Ukraine

\smallskip

 E-mail: rezounenko@univer.kharkov.ua

\end{center}

\begin{quote}
{\bf Abstract.} We investigate a class of non-linear partial
differential equations with discrete state-dependent delays. The
existence and uniqueness of strong solutions for initial functions
from a Banach space are proved. To get the well-posed initial value
problem we restrict our study to a smaller metric space,
construct the dynamical system and prove the existence of a compact global attractor.
\end{quote}

\bigskip


\bigskip


\bigskip
\section{Introduction}\label{sec1}

\marginpar{\tiny Apr 15, 2009}



Theory of dynamical systems is a theory which describes qualitative
properties of systems, changing in time. One of the oldest branches
of this theory is the theory of delay differential equations. We
refer to some classical monographs on the theory of ordinary
(O.D.E.) delay equations \cite{Hale_book,Walther_book,Azbelev}. A
characteristic feature of any type of delay equations is that they
generate infinite dimensional dynamical systems.
 The theory of partial (P.D.E.) delay equations
is essentially less studied since such equations are simultaneously
infinite dimensional in both time (as delay equations) and space (as
P.D.E.s) variables, which makes the analysis more difficult. We
refer to some works which are close to the present research
\cite{travis_webb,Chueshov-JSM-1992,Cras-1995,NA-1998} and to the
monograph \cite{Wu_book}.

Recently, much attention was paid to the investigations of a new
class of delay equations - equations with a state-dependent delays
(SDD) (see e.g.
\cite{Walther2002,Walther_JDE-2003,Walther_JDDE-2007,Krisztin-2003,Nussbaum-Mallet-1992,Nussbaum-Mallet-1996,MalletParet}
and also the survey paper \cite{Hartung-Krisztin-Walther-Wu-2006}
for details and references). The study of these equations
essentially differ from the ones of equations with constant or
time-dependent delays. The main difficulty is that nonlinearities
with SDDs
are not Lipschitz continuous on the space of continuous functions -
the main space of initial data, where the classical theory of delay
equations is developed (see the references above). As a result, the
corresponding initial value problem (IVP), in general, is {\bf not}
well-posed (in the sense of J.~Hadamard
\cite{Hadamard-1902,Hadamard-1932}).
 An explicit example of the non-uniqueness of solutions
to an ordinary equation with state-dependent delay (SDD) is given in
\cite{Winston-1970} (see also
\cite[p.465]{Hartung-Krisztin-Walther-Wu-2006}). As noticed in
\cite[p.465]{Hartung-Krisztin-Walther-Wu-2006} "typically, the IVP
is uniquely solved for initial and other data which satisfy suitable
Lipschitz conditions."

First attempts to study P.D.E.s with SDDs have been made for
different types of delays: for a distributed delay problem in
\cite{Rezounenko-Wu-2006,Rezounenko_JMAA-2007} (see also
\cite{Rezounenko_MMAS-2008})
 and for discrete SDDs in
\cite{Hernandez-2006} (mild solutions, infinite delay) as well as in
\cite{Rezounenko_JMAA-2007} (weak solutions, finite delay).


The following property of solutions of P.D.E.s (with or without
delays) is very important for the study of equations with a discrete
state-dependent delay. Considering any type of solutions (weak,
mild, strong or classical) and having the property $u \in
C([a,b];X),$  one cannot, in general, guarantee that the solution is
a Lipschitz function $u : [a,b] \to X.$ This fact brings essential
difficulties for the extension of the methods developed for O.D.E.s
(see the discussion above).
That is why in previous investigations we proposed alternative
approaches i.e. approximations of a solution of a P.D.E. with a
discrete SDD by a sequence of solutions of P.D.E.s with distributed
SDDs~\cite{Rezounenko-Wu-2006,Rezounenko_JMAA-2007} or use an
"ignoring condition" for a discrete SDD
function~\cite{Rezounenko_NA-2009}.

The main goal of the present work is to 
make a step in extending the approach used for O.D.E.s with SDDs
\cite{Walther2002,Walther_JDE-2003,Hartung-Krisztin-Walther-Wu-2006}
to the case of P.D.E.s. Our idea is to look for a wider space
$Y\supset X$ such that a solution $u : [a,b] \to Y$ be a Lipschitz
function (with respest to the weaker norm of $Y$) and to construct a
dynamical system on a subset of the space  $C([a,b];Y)$. It is
interesting to notice that, in contrast to the previous
investigations, the dynamical system is constructed on a metric
space which is not a linear space.

The article is organized as follows. Section~2 contains formulation
of the model and the proof of the existence and uniqueness of strong
solutions for initial functions from a Banach space. In Section~3 we
construct an evolution operator $S_t$ and study its asymptotic
properties. Here we restrict the evolution operator to a smaller
metric space to get the continuity of $S_t$, which is not available
in the initial Banach space.
 Section~4 deals with the particular case when the delay time is state-independent. Here we also compare the results
with the state-dependent case.


\section{Formulation of the model and basic properties}\label{sec2}



Our  goal is to present an approach to study the following
 partial differential equation with  state-dependent
discrete delay

$$
\frac{\partial }{\partial t}u(t,x)+Au(t,x)+du(t,x)$$
\begin{equation}\label{sdd9-1}=  b\left( [Bu(t-\eta (u_t),\cdot)](x)\right) \equiv \big( F(u_t)
\big)(x),\quad x\in \Omega ,
\end{equation}
 where $A$ is a densely-defined self-adjoint positive linear operator
 with domain $D(A)\subset L^2(\Omega )$ and with compact
  resolvent, so $A: D(A)\to L^2(\Omega )$ generates an analytic semigroup,
  $\Omega $ is a smooth bounded domain in $R^{n_0}$,
 $B: L^2(\Omega ) \to L^2(\Omega ) $ is a bounded operator to be specified later,
$b:R\to R$ is a locally Lipschitz 
map 
and $d$ is a non-negative constant. The function
  $\eta (\cdot): C([-r,0];L^2(\Omega)) \to [0,r]\subset R_{+}$ represents the
state-dependent {\it discrete} delay.
We denote for short $C\equiv C([-r,0];L^2(\Omega)).$ The norms in
$L^2(\Omega)$ and $C$ are denoted by $||\cdot ||$ and $||\cdot ||_C$
respectively. By $\langle \cdot,\cdot\rangle $ we denote the inner
product in $L^2(\Omega)$. As usual for delay equations, we denote by
$u_t$ the function of $\theta\in [-r,0]$ by the formula $u_t\equiv
u_t(\theta)\equiv u(t+\theta).$

\medskip

{\bf Remark~1}. {\it For example, the operator $B$ may be of the
following forms (linear examples)
\begin{equation}\label{sdd9-2} [Bv](x)\equiv  \int_\Omega v(y) \widetilde f(x,y) dy ,\quad x\in
\Omega,
\end{equation}
or even simpler
\begin{equation}\label{sdd9-3} [Bv](x)\equiv
\int_\Omega v(y)  f(x-y) \ell (y) dy ,\quad x\in \Omega,
\end{equation}
where $f : \Omega -\Omega \to R$ is a smooth function, $\ell \in
C^\infty_0 (\Omega)$.
 In the last case the nonlinear term in (\ref{sdd9-1}) takes
the form
\begin{equation}\label{sdd9-4}\big( F(u_t)
\big)(x)\equiv  b\left( \int_\Omega u(t-\eta (u_t), y) f(x-y) \ell (y) dy \right),\quad x\in \Omega ,
\end{equation}
}

We consider equation (\ref{sdd9-1}) with the following initial
conditions
\begin{equation}\label{sdd9-ic}
  u|_{[-r,0]}=\varphi. 
\end{equation}

\medskip
Main assumptions:

{\bf (H.B)} We will need the following Lipschitz property of the
operator $B$
\begin{equation}\label{sdd9-5}
\exists L_B>0 : \forall u,v\in L^2(\Omega ) \Rightarrow ||Bu -
Bv||\le L_B\cdot ||A^{-1/2}(u-v)||.
\end{equation}

{\bf (H.$\eta$)} The discrete delay function $\eta: C\to [0,r]$
satisfies
\begin{equation}\label{sdd9-6}
\exists L_\eta>0 : \forall \varphi,\psi\in C 
\Rightarrow |\eta(\varphi) - \eta(\psi)|\le L_\eta\cdot
\max_{\theta\in [-r,0]} ||A^{-1/2}(\varphi(\theta)-\psi(\theta))||.
\end{equation}

\medskip

{\bf Remark~2}. {\it For the term of the form (\ref{sdd9-3}),
assuming that for all (almost all) $x\in \Omega \Rightarrow f(\cdot
- x)\ell (\cdot)\in D(A^{1/2})$ and $u\in L^2(\Omega)\subset
D(A^{-1/2})$ one gets $|\, \langle u,f(\cdot - x)\ell (\cdot)
\rangle | \le ||A^{-1/2}u||\cdot ||A^{1/2}f(\cdot - x)\ell (\cdot)
||$ which implies
$$ \left( \int_\Omega \big| \int_\Omega u(y)f(y-x)\ell (y) dy \big|^2\, dx \right)^{1/2}
\le ||A^{-1/2}u||\cdot \left( \int_\Omega || A^{1/2} f(\cdot-x)\ell
(\cdot) ||^2\, dx \right)^{1/2}.
$$
Hence, property (H.B) (see (\ref{sdd9-5})) holds with $ L_B\equiv
\left( \int_\Omega || A^{1/2} f(\cdot-x)\ell (\cdot) ||^2\, dx
\right)^{1/2}. $

The same arguments hold (with $ L_B\equiv \left( \int_\Omega ||
A^{1/2} \widetilde f(x,\cdot) ||^2\, dx \right)^{1/2}$) for a more
general term of the form (\ref{sdd9-2}).
}
\medskip


%
%

Now we introduce the following

\medskip

{\bf Definition~1.} {\it A vector-function $u(t)\in
C([-r,T];D(A^{-1/2}))\cap C([0,T];D(A^{1/2})) \cap L^2(0,T;D(A))$
with derivative $\dot
u(t)\in L^\infty (0,T;D(A^{-1/2}))$ 
is a strong solution of problem (\ref{sdd9-1}), (\ref{sdd9-ic}) on
an interval $[0,T]$ if
\begin{itemize}
    \item $u(\theta)=\varphi (\theta)$ for $\theta\in [-r,0]$;
    \item for any function $v\in L^2(0,T;L^2(\Omega)) $ such that $\dot v\in
L^2(0,T;D(A^{-1})) $ and $v(T)=0,$ one has
\begin{equation}\label{sdd9-7}
\hskip-10mm -\int^T_0 \langle u(t),\dot v(t)\rangle \, dt  +
\int^T_0 \langle A^{1/2} u(t), A^{1/2} v(t)\rangle \, dt= \langle
\varphi (0),v(0)\rangle + \int^T_0 \langle F(u_t)-du(t), v(t)\rangle
\, dt.
\end{equation}
\end{itemize}
}%

\medskip

Let us introduce the following space 
\begin{equation}\label{sdd9-8}
{\mathcal L}\equiv \left\{ \varphi \in
C([-r,0];D(A^{-1/2}))\, |\, \sup\limits_{s\neq t} \left\{
\frac{||A^{-1/2}(\varphi (s)-\varphi (t)) ||}{|s-t|} \right\}
<+\infty;\, \varphi(0)\in D(A^{1/2})\right\}
\end{equation}

with the natural norm
\begin{equation}\label{sdd9-9}
\hskip-10mm ||\varphi||_{\mathcal L}\equiv \max_{s\in [-r,0]}
||A^{-1/2}\varphi(s)|| +\sup\limits_{s\neq t} \left\{
\frac{||A^{-1/2}(\varphi (s)-\varphi
 (t)) ||}{|s-t|} \right\} + ||A^{1/2} \varphi(0)||.
\end{equation}

\medskip

Now we prove the following theorem on the existence and uniqueness of  
solutions.
\medskip

{\bf Theorem~1.} {\it Let assumptions (H.B) and (H.$\eta$) hold (see
(\ref{sdd9-5}), (\ref{sdd9-6})). Assume the function $b:R\to R$ is
locally Lipschitz and bounded ($b(\cdot)\le M_b$).


 Then for any initial function $\varphi \in {\mathcal L}$ (the space ${\mathcal L}$ is defined in (\ref{sdd9-8})) the problem (\ref{sdd9-1}), (\ref{sdd9-ic}) has a unique strong
solution on any time interval $[0,T]$. The solution has the property
$\dot u \in L^2 (0,T;L^2 (\Omega)).$
}%
\medskip

{\bf Remark~3.} {\it Let us notice that we do not assume that
$\varphi\in L^2 (-r,0;D(A))$, but the definition of strong solution
above implies that 
\begin{equation}\label{sdd9-11}
u_t\in L^2 (-r,0;D(A)), \quad \forall t\ge r.
\end{equation}
}%

\medskip
{\it Proof of theorem~1.} Let us denote by $\{ e_k\}^\infty_{k=1}$
an orthonormal basis of $L^2(\Omega)$ such that $Ae_k=\lambda_ke_k$,
$0< \lambda_1<\ldots<\lambda_k\to +\infty$.

Consider Galerkin approximate solutions of order $m$ :\\
$u^m=u^m(t,x)=\sum^m_{k=1} g_{k,m}(t) e_k,$
such that 
 \begin{equation}\label{sdd9-12}
\left\{ \begin{array}{ll} &\langle \dot u^m+Au^m +du^m-F(u^m_{t}),
e_k\rangle =0,\\
&
\langle u^m(\theta),e_k\rangle=\langle \varphi (\theta) , e_k\rangle
,\,\,\forall \theta\in [-r,0] \end{array}\right.
 \end{equation}
$\forall k=1,\ldots,m$. Here $g_{k,m}\in C^1(0,T;R)\cap L^2(-r,T;R)$
with $\dot g_{k,m}(t)$ being absolutely continuous.

The system (\ref{sdd9-12}), is an (ordinary) differential equation
in $R^m$ with a concentrated state-dependent delay for the unknown
vector function $U(t)\equiv (g_{1,m}(t), \ldots,g_{m,m}(t))$ (the
corresponding theory is developed in
\cite{Walther_JDE-2003,Walther_JDDE-2007} see also a recent review
\cite{Hartung-Krisztin-Walther-Wu-2006}).

The key difference between equations with state-dependent and
state-independent delays is that the first type of equations is not
well-posed in the space of continuous (initial) functions. To get
the well-posed initial value problem, the theory
\cite{Walther_JDE-2003,Walther_JDDE-2007,Hartung-Krisztin-Walther-Wu-2006}
suggests to restrict considerations to a smaller space of Lipschitz
continuous functions or even to a smaller subspace of
$C^1([-r,0];R^m).$

Condition $\varphi \in {\mathcal L}$ implies that initial data
$U(\cdot)|_{[-r,0]}\equiv P_m \varphi (\cdot) $ is Lipschitz
continuous as a function from $[-r,0]$ to $R^m.$ Here $P_m$ is the
orthogonal projection onto the subspace $span\, \{ e_1,\ldots,
e_m\}\subset L^2(\Omega).$  Hence we can apply the theory of O.D.E.s
with state-dependent delay (see e.g.
\cite{Hartung-Krisztin-Walther-Wu-2006}) to get the local existence
and uniqueness of solutions of (\ref{sdd9-12}).

Now we look for an a priory estimate to prove the continuation of
solutions $u^m$
of (\ref{sdd9-12}) on any time interval $[0,T]$ and then use it for the proof (by the method of compactness, see \cite{Lions}) of the existence of strong solutions to (\ref{sdd9-1}), (\ref{sdd9-ic}).  

We multiply the first equation in (\ref{sdd9-12}) by $\lambda_k
g_{k,m}$ and sum for $k=1,\ldots,m$ to get
$${1\over 2} {d \over dt} ||A^{1/2}u^m(t)||^2 + ||Au^m(t)||^2 + d\cdot ||A^{1/2}u^m(t)||^2= \langle P_m F(u^m_t),Au^m(t)\rangle
$$
$$\le {1\over 2} ||P_m F(u^m_t)||^2 + {1\over 2} ||Au^m(t)||^2.
$$
Since the function $b$ is bounded ($b(\cdot)\le M_b$), we have
$||F(u^m_t)||^2\le M^2_b |\Omega|$ (here $|\Omega|\equiv \int_\Omega
1 \, dx$) and, as a result, we conclude that
\begin{equation}\label{sdd9-13}
{d \over dt} ||A^{1/2}u^m(t)||^2 + ||Au^m(t)||^2 \le M^2_b |\Omega|.
\end{equation}
We integrate (\ref{sdd9-13}) with respect to $t$, use the properties
$\varphi (0)\in D(A^{1/2}),$ 
$u^m(0)=P_m\varphi (0)\in D(A^{1/2})$ and
$||A^{1/2}u^m(0)||=||A^{1/2}P_m\varphi (0)|| \le ||A^{1/2}\varphi
(0)||$ to get an a priory estimate
\begin{equation}\label{sdd9-14}
\hskip-10mm  ||A^{1/2}u^m(t)||^2 + \int^t_0 ||Au^m(\tau)||^2 \,
d\tau \le ||A^{1/2}\varphi (0)||^2 + M^2_b |\Omega|\cdot T, \quad
\forall m, \forall t\in [0,T].
\end{equation}

Estimate (\ref{sdd9-14}) means that
$$ \{ u^m \}^\infty_{m=1} \hbox{ is a bounded set in }  L^\infty(0,T;D(A^{1/2}))\cap L^2(0,T;D(A)).
$$
Using this and (\ref{sdd9-12}), we get that
$$ \{ \dot u^m \}^\infty_{m=1} \hbox{ is a bounded set in }  L^\infty(0,T;D(A^{-1/2}))\cap L^2(0,T;L^2(\Omega)).
$$
Hence the family $\{ (u^m;\dot u^m ) \}^\infty_{m=1}$ is a bounded
set in
\begin{equation}\label{sdd9-15}
\hskip-15mm Z_1\equiv \left( L^\infty(0,T;D(A^{1/2}))\cap
L^2(0,T;D(A))\right) \times \left( L^\infty(0,T;D(A^{-1/2}))\cap
L^2(0,T;L^2(\Omega)) \right).
\end{equation}
Therefore there exist a subsequence $\{ (u^k;\dot u^k ) \}$ and an
element  $(u;\dot u )\in Z_1$ such that
\begin{equation}\label{sdd9-16}
\{ (u^k;\dot u^k ) \} \hbox{ *-weak converges to } (u;\dot u )
\hbox{ in } Z_1.
\end{equation}
The proof that any *-weak limit is a strong solution is standard.

\medskip

Now we prove the uniqueness of strong solutions.

Using the properties
$\varphi \in {\mathcal L}$,  the definition of a strong solution $v$
and $\dot v(t) \in L^\infty (0,T; D(A^{-1/2}))$ (see
(\ref{sdd9-16})) we have that for any such a solution $v$ and any
$T>0$ there exists $L_{v,T}>0,$ such that
\begin{equation}\label{sdd9-18}
||A^{-1/2}(v(s^1)-v(s^2))||\le L_{v,T} \cdot |s^1-s^2|, \quad
\forall s^1,s^2 \in [-r,T].
\end{equation}

Consider two strong solutions $u$ and $v$ of (\ref{sdd9-1}),
(\ref{sdd9-ic}) (not necessarily with the same initial function).

Assumption (H.B) (see (\ref{sdd9-5})) and the Lipschitz property of
$b$ imply
$$ || F(u_{s^1})- F(v_{s^2}) ||^2=\int_\Omega | \, b\left( [Bu](s^1-\eta(u_{s^1}),x)\right)
- b\left( [Bv](s^2-\eta(u_{s^2}),x)\right) |^2 \, dx $$
$$\le L^2_b \int_\Omega | \,[B u](s^1-\eta(u_{s^1}),x) - [B v] (s^2-\eta(u_{s^2}),x) |^2 \, dx
$$
$$= L^2_b\cdot || \,[B u](s^1-\eta(u_{s^1}),\cdot) - [B v]
(s^2-\eta(u_{s^2}),\cdot) ||^2
$$
\begin{equation}\label{sdd9-17}
\le L^2_b L^2_B\cdot
||A^{-1/2}(u(s^1-\eta(u_{s^1}))-v(s^2-\eta(u_{s^2}))||^2.
\end{equation}

Now, for any two strong solutions, we have
$$  F(u_t)- F(v_t)  = b(Bu(t-\eta (u_t))) - b(Bv(t-\eta (v_t))) \pm b(Bv(t-\eta (u_t))).
$$
Using the Lipschitz properties of $b, B$ and $\eta$ (see
(\ref{sdd9-5}), (\ref{sdd9-6})), and also (\ref{sdd9-18}),
(\ref{sdd9-17}), one gets

$$\hskip-110mm ||F(u_t)- F(v_t) || $$
$$ \le L_b L_B \left( \max_{s\in [t-r,t]} ||A^{-1/2}(u(s)-v(s))|| +
||A^{-1/2}(v(t-\eta (u_t))-v(t-\eta (v_t)))|| \right)
$$
$$\le L_b L_B \left(  ||A^{-1/2}(u_t-v_t)||_C +
L_{v,T} \cdot |\eta (u_t))-\eta (v_t)| \right)
$$
\begin{equation}\label{sdd9-19}
 \le
L_b L_B \left(  1 + L_{v,T} \cdot L_\eta \right) \cdot
||A^{-1/2}(u_t-v_t)||_C.
\end{equation}
We denote for short
\begin{equation}\label{sdd9-19a} C_{v,T}\equiv L_b L_B \left(
1 + L_{v,T} \cdot L_\eta \right).
\end{equation}

Now the standard variation-of-constants formula\\ $ u(t)=e^{-At}
u(0) + \int^t_0 e^{-A(t-\tau)} F(u_\tau)\, d\tau$ and
(\ref{sdd9-19}) give
$$
||A^{-1/2}(u_t - v_t)||_C \le ||A^{-1/2}(u_0 - v_0)||_C +  C_{v,T}
\cdot \int^t_0 e^{-\lambda_1 (t-\tau)} ||A^{-1/2}(u_\tau -
v_\tau)||_C \, d\tau.
$$
The last estimate (by Gronwall's lemma) implies
\begin{equation}\label{sdd9-21}
\hskip-10mm ||A^{-1/2}(u_t - v_t)||_C \le ||A^{-1/2}(u_0 - v_0)||_C
\cdot \left[ 1+ \frac{C_{v,T}}{C_{v,T}-\lambda_1}
\left(e^{(C_{v,T}-\lambda_1) t} -1\right) \right],
\end{equation}
which gives the uniqueness of strong solutions of (\ref{sdd9-1}),
(\ref{sdd9-ic}).

The proof of theorem~1 is complete. \finproof

\medskip

{\bf Remark~4.} {\it It is very important that the term $\left[ 1+
\frac{C_{v,T}}{C_{v,T}-\lambda_1} \left(e^{(C_{v,T}-\lambda_1) t}
-1\right) \right] $ in (\ref{sdd9-21}}) tends to $+\infty$, when
$L_{v,T}\to +\infty$, except the case  $L_\eta = 0$ (see
(\ref{sdd9-19a})).

\medskip

Let us get an additional estimate for strong solutions.

The standard variation-of-constants formula $ u(t)=e^{-At} u(0) +
\int^t_0 e^{-A(t-\tau)} F(u_\tau)\, d\tau$, (\ref{sdd9-19}),
(\ref{sdd9-19a}) and the estimate  $|| A^\alpha e^{-tA}|| \le
\left({\alpha\over t}\right)^\alpha e^{-\alpha}$ (see e.g.
\cite[(1.17), p.84]{Chueshov_book}) give
$$||A^{1/2}(u(t) - v(t))|| \le e^{-\lambda_1t}||A^{1/2}(u(0) - v(0))|| +
\int^t_0 || A^{1/2} e^{-A(t-\tau)}|| \cdot ||F(u_\tau) -F(v_\tau) ||
\, d\tau $$
 \begin{equation}\label{sdd9-20}
 \le
e^{-\lambda_1t}||A^{1/2}(u(0) - v(0))|| + 2 t^{1/2}\left({1\over
2}\right)^{1/2} e^{-1/2}\cdot C_{v,T} \cdot ||A^{-1/2}(u_0 -
v_0)||_C.
\end{equation}
Here we used $|| A^{1/2} e^{-A(t-\tau)}|| \le \left({1/2\over
t-\tau}\right)^{1/2} e^{-1/2}$ and $\int^t_0 (t-\tau)^{-1/2}d\tau =
2 t^{1/2}.$

Now estimates (\ref{sdd9-21}), (\ref{sdd9-20}) give

$$||A^{1/2}(u(t) - v(t))|| + ||A^{-1/2}(u_t - v_t)||_C $$
 \begin{equation}\label{sdd9-20a}
 \le
e^{-\lambda_1t}||A^{1/2}(u(0) - v(0))|| + D_{v,T} \cdot
||A^{-1/2}(u_0 - v_0)||_C.
\end{equation}
Here we denote
 \begin{equation}\label{sdd9-31} D_{v,T} \equiv 2 T^{1/2}\left({1\over
2}\right)^{1/2} e^{-1/2}\cdot C_{v,T} + \left[ 1+
\frac{C_{v,T}}{C_{v,T}-\lambda_1} \left(e^{(C_{v,T}-\lambda_1) T}
-1\right) \right].
\end{equation}

%
\section{Asymptotic behavior}\label{sec3}
%

In this section we study long-time behavior of the strong solutions
of the problem (\ref{sdd9-1}), (\ref{sdd9-ic}).

Due to theorem~1, we define in the standard way the evolution
semigroup $S_t : {\mathcal L} \to {\mathcal L}$ (the space
${\mathcal L}$ is defined in (\ref{sdd9-8})) by the formula
\begin{equation}\label{sdd9-22}
S_t \varphi \equiv u_t,\quad t\ge 0,
\end{equation}
where $u(t)$ is the unique strong solution of the problem
(\ref{sdd9-1}), (\ref{sdd9-ic}).

\medskip

 {\bf Remark~5}. {\it We emphasize that the evolution semigroup $S_t : {\mathcal L} \to {\mathcal L}$ is not a
dynamical system in the standard sense (see
e.g.\cite{Babin-Vishik,Temam_book,Chueshov_book}) since $S_t$ is not
a continuous mapping in the topology of ${\mathcal L}$ i.e. the
problem (\ref{sdd9-1}), (\ref{sdd9-ic}) is not well-posed in the
sense of J.~Hadamard~\cite{Hadamard-1902,Hadamard-1932}.
}%

\medskip

Our first goal is to prove

\medskip

{\bf Lemma~1.}  {\it Let all the assumptions of theorem~1 be
satisfied. Then for any $\alpha \in ({1\over 2},1)$, there exists a
bounded in the space $C^1([-r,0];D(A^{-1/2}))\cap
C([-r,0];D(A^{\alpha}))$ set ${\mathcal B}V_\alpha$, which absorbs
any strong solution of the problem (\ref{sdd9-1}), (\ref{sdd9-ic})
with any initial function $\varphi \in {\mathcal L}.$
}%

\medskip {\it Proof of lemma~1.}
Using $||A^{1/2}v||^2\le \lambda^{-1}_1\cdot ||Av||^2,$ we get from
(\ref{sdd9-13}) that
$${d \over dt} ||A^{1/2}u^m(t)||^2 + \lambda_1 ||A^{1/2}u^m(t)||^2 \le M^2_b |\Omega|.$$
We multiply the last estimate by $e^{\lambda_1 t}$ and integrate
over $[0,t]$ to obtain
 $$||A^{1/2}u^m(t)||^2\le ||A^{1/2}u^m(0)||^2 e^{-\lambda_1 t} +
\lambda^{-1}_1M_b^2 |\Omega|
$$
\begin{equation}\label{sdd9-30}\le ||A^{1/2}\varphi(0)||^2
e^{-\lambda_1 t} + \lambda^{-1}_1M_b^2 |\Omega|.
\end{equation}
 This and
(\ref{sdd9-12}) give $||A^{-1/2}\dot u^m(t)||^2 \le
2||A^{1/2}\varphi(0)||^2 e^{-\lambda_1 t} + 2\lambda^{-1}_1M_b^2
|\Omega|+ M_b^2 |\Omega|.$ The last two estimates imply
\begin{equation}\label{sdd9-23}
 \hskip-10mm   ||A^{1/2}u^m(t)||^2+||A^{-1/2}\dot u^m(t)||^2\le 3 ||A^{1/2}\varphi(0)||^2 e^{-\lambda_1 t} +
(1+3\lambda^{-1}_1) M_b^2 |\Omega|.
\end{equation}


 We get an analogous estimate for a strong solution of the
problem (\ref{sdd9-1}), (\ref{sdd9-ic}), using the well-known

\smallskip

\noindent {\bf Proposition~1.} \cite[Theorem~9]{yosida}.
 {\it Let $X$ be a Banach space. Then any *-weak convergent sequence
 $\{ w_k\}^\infty_{n=1}\in X^{*}$  *-weak converges to an element
 $w_\infty\in X^{*}$ and $\Vert w_\infty\Vert_X \le\liminf_{n\to\infty} \Vert w_n\Vert_X.$
}

\medskip

\noindent


Now we consider the space $V\equiv C^1([-r,0];D(A^{-1/2}))\cap
C([-r,0];D(A^{1/2}))$, fix any positive $\varepsilon_0$ and obtain
that the ball ${\mathcal B}_0$ of $V$
\begin{equation}\label{sdd9-24}
  {\mathcal B}_0\equiv \left\{  v\in V:  ||v||^2_V \le R^2_0\equiv
(1+3\lambda^{-1}_1) M_b^2 |\Omega| + \varepsilon_0 \right\}
\end{equation}
is absorbing for any strong solution of the problem (\ref{sdd9-1}),
(\ref{sdd9-ic}) (see (\ref{sdd9-23})).

Now we are in a position to use the arguments presented in
\cite[Lemma 2.4.1, p.101]{Chueshov_book} and get (for any ${1\over
2}<\alpha <1$) the existence of the absorbing ball
\begin{equation}\label{sdd9-27}
 {\mathcal B}_\alpha \equiv \left\{  v\in C([-r,0];D(A^{\alpha})) :  ||v||_{C([-r,0];D(A^{\alpha}))}
 \le R_\alpha
 \right\},
\end{equation}
where $R_\alpha \equiv (\alpha -1/2)^{\alpha -1/2}\cdot
\left[\lambda^{-1/2}_1 M_b \sqrt{|\Omega|}+\varepsilon\right] +
{\alpha^\alpha \over 1-\alpha} \cdot M_b \sqrt{|\Omega|}$ with any
fixed $\varepsilon>0.$ More precisely, the standard
variation-of-constants formula $ u(t)=e^{-At} u(0) + \int^t_0
e^{-A(t-\tau)} F(u_\tau)\, d\tau$ and the estimate  $|| A^\alpha
e^{-tA}|| \le \left({\alpha\over t}\right)^\alpha e^{-\alpha}$ (see
e.g. \cite[(1.17), p.84]{Chueshov_book}) give
$$ ||A^{\alpha} u(t+1)|| \le (\alpha -1/2)^{\alpha -1/2} ||A^{1/2}
u(t)|| + \int^{t+1}_t \left( \alpha \over t+1-\tau \right)^\alpha
||F(u_\tau)||\, d\tau.
$$
Let us consider any bounded in ${\mathcal L}$ set $\hat B$.
 Estimate (\ref{sdd9-30}) and
Proposition~1 give $||A^{1/2} u(t)||\le \left[\lambda^{-1/2}_1 M_b
\sqrt{|\Omega|}+\varepsilon\right]$ for all $t\ge t_{\hat B}$ (here
$t_{\hat B}$ depends on ${\hat B}$ only). These and the estimate
$||F(u_\tau)|| \le M_b \sqrt{|\Omega|}$ imply (\ref{sdd9-27}).

The above estimates (\ref{sdd9-24}), (\ref{sdd9-27})  show that
there exists a subset (a ball) of $V_\alpha \equiv
C^1([-r,0];D(A^{-1/2}))\cap C([-r,0];D(A^{\alpha}))$ (here ${1\over
2}<\alpha <1$)
\begin{equation}\label{sdd9-28}
 {\mathcal B}V_\alpha \equiv \left\{  v\in V_\alpha :  ||v||_{V_\alpha}
 \le \widehat R_\alpha
 \right\},
\end{equation}
such that for any strong solution, starting in $\varphi$ from any
bounded set ${\hat B}\subset {\mathcal L},$ there exists $t_{\hat
B}\ge 0$ such that
\begin{equation}\label{sdd9-29}
S_t \varphi \in {\mathcal B}V_\alpha, \qquad \hbox{ for all } \qquad
t\ge t_{\hat B}.
\end{equation}

The proof of lemma~1 is complete. \rule{5pt}{5pt}


\medskip

We will use notation

$$ |||\varphi||| \equiv \sup\limits_{s\neq t} \left\{ \frac{||A^{-1/2}(\varphi (s)-\varphi
(t)) ||}{|s-t|} \right\} \hbox{ for } \varphi \in  {\mathcal L}.
$$

Let us fix $R^0>0$ and consider the metric space ${\mathcal
L}_{R^0}$
%
%
which is the set \\ $\left\{ \varphi \in  {\mathcal L} \, : \,
|||\varphi||| \le R^0\right\} $ equipped with the metrics (c.f.
(\ref{sdd9-9}))
\begin{equation}\label{sdd9-26}
\rho (\varphi,\phi) \equiv \max_{s\in [-r,0]}
||A^{-1/2}(\varphi(s)-\phi(s))||  + ||A^{1/2}
(\varphi(0)-\phi(0))||.
\end{equation}

One can check that $({\mathcal L}_{R^0};\rho)$ is a complete metric
space and any
%
set \\ $\left\{ \varphi \in  {\mathcal L} \, : \, |||\varphi||| \le
R^1< R^0\right\} $ is closed.

\medskip

We need the following (technical) assumption

\medskip

{\bf (H.I)} {\it There exists $R^0> \widehat R_\alpha $ ($\widehat
R_\alpha $ is defined in (\ref{sdd9-28})) such that the set

$\left\{ \varphi \in {\mathcal L} \, : \, |||\varphi||| \le
R^0\right\} $ is positively invariant for the semigroup $S_t$
i.e.}
\begin{equation}\label{sdd9-32}
\forall \varphi\in {\mathcal L} : |||\varphi|||  \le R^0 \quad
\Rightarrow \quad  |||S_t \varphi||| \le R^0,\quad \forall t>0.
\end{equation}

\bigskip

Our next result is the following

\bigskip

 {\bf Theorem~2.}  {\it Let (H.I) and all the assumptions of theorem~1 be
 satisfied. Then the evolution
semigroup $S_t : {\mathcal L}_{R^0} \to {\mathcal L}_{R^0}$, (see
(\ref{sdd9-22})),
possesses a global attractor in the metric space $({\mathcal L}_{R^0};\rho).$ 
 }%

\medskip

{\it Proof of theorem~2.}
Now we concentrate on the metric space $({\mathcal L}_{R^0};\rho)$
(here $R^0>\widehat R_\alpha$). The reason for this is that the
evolution semigroup $S_t$ is not continuous on the whole space
${\mathcal L}$ (see remark~5). On the other hand, we notice:

\medskip

 {\bf Remark~6}. {\it Estimate (\ref{sdd9-20a}) implies that
 the evolution semigroup  $S_t$ is a continuous mapping
in the topology of $({\mathcal L}_{R^0};\rho)$ i.e. $\rho
(S_t\varphi, S_t\phi) \le D_{v,T} \cdot \rho (\varphi, \phi)$ for
$\varphi,\phi \in {\mathcal L}_{R^0}$, and $t\in [0,T].$ Here
$D_{v,T}$ is defined by (\ref{sdd9-31}) (see also (\ref{sdd9-19a}))
with $L_{v,T}=R^0$  (c.f. (\ref{sdd9-18})).
}%

\medskip


Corollary~4 from \cite{Simon} implies that ${\mathcal B}V_\alpha $
is relatively compact in $C([-r,0];D(A^{-1/2}))$ (see also
\cite[lemma~1]{Simon}). This fact and the property $||A^\alpha
\varphi (0)|| \le \widehat R_\alpha, {1\over 2} < \alpha <1$ for all
$\varphi\in {\mathcal B}V_\alpha $ gives that ${\mathcal B}V_\alpha
$ is relatively compact in the topology of $({\mathcal
L}_{R^0};\rho)$.

Let us consider the following set
$$ K\equiv Cl\, [{\mathcal B}V_\alpha]_{({\mathcal L}_{R^0};\rho)},$$
where $Cl\, [\cdot ]_{({\mathcal L}_{R^0};\rho)}$ is the closure in
the topology of $({\mathcal L}_{R^0};\rho)$. The above properties
show that $K$ is compact in $({\mathcal L}_{R^0};\rho)$.

We get (see (\ref{sdd9-29})) that for any strong solution, starting
in $\varphi$ from any bounded set ${\widetilde B}\subset {\mathcal
L}_{R^0}$, there exists $t_{\widetilde B}\ge 0$ such that $$ S_t
\varphi \in {\mathcal B}V_\alpha \subset K, \qquad \hbox{ for all }
\qquad t\ge t_{\widetilde B}.$$

As a result, we conclude that the evolution semigroup $S_t$ is
asymptotically compact (and dissipative) in $({\mathcal
L}_{R^0};\rho)$.

 Finally, by the classical theorem on the existence of an
attractor (see, for example,
\cite{Babin-Vishik,Temam_book,Chueshov_book}) one gets that
$(S_t;({\mathcal L}_{R^0};\rho))$ has a compact global attractor.
The proof of Theorem~2 is complete. \rule{5pt}{5pt}

\subsection{Equation with a modified nonlinearity}


Discussing the technical assumption (H.I), we notice that even in
the case when (H.I) is not satisfied for the original system,
Lemma~1 allows one to consider a modified system without modifying
the long term dynamics of $S_t$ (see
\cite{Debussche-Temam_DCDS-1996}). More precisely, one chooses
\cite[p.545]{Debussche-Temam_DCDS-1996} a $C^\infty$ function $\chi
: [0,+\infty) \to [0,1]$ such that
$$\left\{\begin{array}{ll}
    \chi (s)=1, & s\in [0,1]; \\
    \chi (s)=0, & s\in [2,+\infty); \\
    0\le \chi (s)\le 1, & s\in [1,2] \\
  \end{array}
  \right.
$$
and set
$$ \widetilde F(\varphi)\equiv \chi \left( { ||\varphi||_{H,d} \over \widehat
R_\alpha}\right)\cdot F(\varphi).
$$
Here we denoted for short $||\varphi||_{H,d} \equiv ||A^{1/2}\varphi
(0)|| + d\cdot ||A^{-1/2}\varphi ||_C.$

As a result, the modified system (\ref{sdd9-1}) (with $ \widetilde
F(\varphi)$ instead of $ F(\varphi) $) has the same behavior inside
of the (absorbing) set ${\mathcal B}V_\alpha $ (in fact, the
behavior is unchanged inside of a bigger set $\{ \varphi :
||\varphi||_{H,d} \le  \widehat R_\alpha \}$).

Now, using $||\widetilde F(\varphi)|| \le ||F(\varphi)|| \le M_b
\sqrt{|\Omega|}$ and the estimate $||A^{-1/2}\dot u(t) || \le
||A^{1/2} u(t) || + d ||A^{-1/2} u(t) || + M_b \sqrt{|\Omega|}$, we
conclude that the set
$$ {\mathcal  L}(\widehat R_\alpha ) \equiv \left\{
\varphi \in {\mathcal  L}\quad : \quad ||\varphi||_{H,d} \le
2\widehat R_\alpha;\quad |||\varphi |||\le 2\widehat R_\alpha + M_b
\sqrt{|\Omega|}\,\, \right\}
$$
is positively invariant for the evolution operator, constructed by
solutions of (\ref{sdd9-1}) with the modified nonlinearity $
\widetilde F(\varphi)$. We notice that ${\mathcal B}V_\alpha \subset
{\mathcal  L}(\widehat R_\alpha ) \subset \{ \varphi :
||\varphi||_{H,d} \le  2\widehat R_\alpha \}.$

This invariantness of the set $ {\mathcal  L}(\widehat R_\alpha )$
gives the possibility to define (similar to (\ref{sdd9-22})) an
evolution operator $\widetilde S_t : {\mathcal  L}(\widehat R_\alpha
) \to {\mathcal  L}(\widehat R_\alpha ) $ by solutions of
(\ref{sdd9-1}) with the modified nonlinearity $\widetilde F.$

Following the line of arguments presented in theorem~2, we prove the
following analog to theorem~2

\medskip

 {\bf Theorem~3.}  {\it Let (H.I) and all the assumptions of theorem~1 be
 satisfied. Then the evolution
semigroup $\widetilde S_t : {\mathcal  L}(\widehat R_\alpha ) \to
{\mathcal  L}(\widehat R_\alpha ) $, possesses a global attractor in
the metric space $({\mathcal  L}(\widehat R_\alpha );\rho).$ Here,
as before, $\rho$ is the metrics defined by (\ref{sdd9-26}).
 }%

\medskip

{\bf Remark~7}. {\it Discussing the restriction of our study from
the linear space ${\mathcal L}$ to the metric spaces $({\mathcal
L}_{R^0};\rho)$ or $({\mathcal  L}(\widehat R_\alpha );\rho)$, we
notice that it is a natural step even for ordinary differential
equations with a discrete state-dependent delay. For example, in
\cite[Proposition~1 and Corollary~1]{Walther2002} 
it is shown that maximal solutions of a scalar delay equation with a
SDD constitute a semiflow on the set $\{ \phi : \hbox{Lip}(\phi)\le
k, || \phi ||<w \}\subset C([-r,0],R).$ Here Lip$(\phi)=\sup_{x\neq
y} |\phi(x)-\phi(y)|\cdot |x-y|^{-1}.$
} 


\section{Particular case of a state-independent delay
($\eta=\hbox{const}$)}\label{sec4}
%

In this particular case, the assumption (H.$\eta$) (see
(\ref{sdd9-6})) is valid automatically with $L_\eta=0.$ Following
the proof of theorem~1, one can see that the assumption
$\sup\limits_{s\neq t} \left\{ ||A^{-1/2}(\varphi (s)-\varphi (t))
||\cdot |s-t|^{-1} \right\}< +\infty $ is not needed in the case
$\eta=\hbox{const}$. This implies that for any initial  function
$\varphi \in H$ (c.f. (\ref{sdd9-8})),
\begin{equation}\label{sdd9-33}
  H\equiv \left\{ \varphi \in C([-r,0];D(A^{-1/2})) \, | \, \, \varphi(0)\in
D(A^{1/2}) \right\}
\end{equation}
the problem (\ref{sdd9-1}), (\ref{sdd9-ic}) has a strong solution.
The uniqueness of a strong solution follows from (\ref{sdd9-20a})
and the fact that $L_\eta=0$ implies $D_{v,T}$ (defined in
(\ref{sdd9-31})) is bounded for any $\varphi \in H$ (c.f. remark~4
and (\ref{sdd9-19a})). This fact gives the continuity of $S_t : H\to
H$ (c.f. remark~5) and as a consequence, that the pair $(S_t;H)$ is
a dynamical system.

Following the proofs of lemma~1 and theorem~2 we have the following
result.

\medskip

{\bf Theorem~4.} {\it Assume $\eta=\hbox{const}$. Let the assumption
(H.B)  hold and the function $b:R\to R$ be locally Lipschitz and
bounded.


 Then for any initial function $\varphi \in H$ the problem (\ref{sdd9-1}), (\ref{sdd9-ic}) has a unique strong
solution on any time interval $[0,T]$. The solution has the property
$\dot u \in L^2 (0,T;L^2 (\Omega)).$

Moreover, the pair $(S_t;H)$ constitutes  a dynamical system which
possesses a global attractor. The attractor is a bounded set in
$C^1([-r,0];D(A^{-1/2}))\cap C([-r,0];D(A^{\alpha}))$ for any
$\alpha \in ({1\over 2},1)$.
}%
\medskip

Now we can compare two cases (state-dependent and state-independent
delays), assuming
\begin{itemize}
  \item the assumption (H.B) holds and
  \item the function $b:R\to
R$ is locally Lipschitz and bounded.
\end{itemize}

\bigskip

\noindent\begin{tabular}{|c|c|c|}
  \hline
   & state-{\bf dependent} delay $\eta$ & state-{\bf independent} delay  \\
  \hline
  The existence  &  &  \\
  and uniqueness of solutions & $\varphi\in {\mathcal L}\subset H$ and (H.$\eta$) & $\varphi\in H$ \\
   &  &  \\
  \hline
    &  &  \\
  The continuity of $S_t$ and & $S_t : ({\mathcal L}_{R^0};\rho) \to ({\mathcal L}_{R^0};\rho) $  &
    \\
  existence of an attractor & $\widetilde S_t :  ({\mathcal  L}(\widehat R_\alpha );\rho) \to ({\mathcal  L}(\widehat R_\alpha );\rho)$ & $S_t : H\to H$ \\
  \hline
\end{tabular}

\bigskip


{\bf Remark~8}. {\it We notice that ${\mathcal L}_{R^0}\subset
{\mathcal L}\subset H$ and the metrics $\rho$ is the natural metrics
of the space~$H.$}

\bigskip

As an application (for both cases of state-dependent and
state-independent delays) we can consider the diffusive Nicholson's
blowflies equation (see e.g. \cite{So-Yang,So-Wu-Yang}) with
state-dependent delays. More precisely, we consider equation
(\ref{sdd9-1}) where $-A$ is the Laplace operator with the Dirichlet
boundary conditions, $\Omega\subset R^{n_0}$ is a bounded domain
with a smooth boundary, the function $f$ can be,
 for example, $ f(s)={1\over \sqrt{4\pi\alpha}}
e^{-s^2/4\alpha}$, as in \cite{So-Wu-Zou} (see remark~2),
the nonlinear (birth) function $b$ is
given by $b(w)=p\cdot we^{-w}.$ Function $b$ is bounded
, so for any delay function $\eta$, satisfying $(H.\eta)$, the
conditions of Theorems~1,2 are valid (in the case when the
assumption (H.I) is satisfied). As a result, we conclude that
the initial value problem (\ref{sdd9-1}),(\ref{sdd9-ic}) is
well-posed in $({\mathcal L}_{R^0};\rho)$ and the dynamical system
$(S_t,({\mathcal L}_{R^0};\rho))$ has a global attractor
(Theorem~2).

If necessary, we modify the system according to subsection 3.1 and
get the existence of a global attractor for the dynamical system
$\left( \widetilde S_t; ({\mathcal L}(\widehat R_\alpha );\rho)
\right)$.

\bigskip



\medskip

\noindent {\bf Acknowledgements.}  The author wishes to thank
 I.D.~Chueshov for useful discussions of an early version of the
manuscript.

\vskip2cm




\hfill Kharkiv

\bigskip

\hfill April 15, 2009


\end{document}